\documentclass[11pt]{article}

\input{xypic.tex} \xyoption{all}
\usepackage{amsthm,amssymb,latexsym,amsmath,amscd}
\usepackage{graphics}

\newcommand{\C}{\mathbb{C}}
\newcommand{\Cl}{\mathcal{C} \ell}
\newcommand{\Ccl}{\C \ell}

\newcommand{\E}{\mbox{\bf E}}
\newcommand{\grad}{\mbox{grad}}

\newcommand{\Hh}{\mbox{\bf H}}

\newcommand{\mb}{\mathbb}
\newcommand{\mc}{\mathcal}
\newcommand{\id}{{\mathbf 1}}
\newcommand{\R}{\mathbb{R}}

\newcommand{\Sc}{\mathbb{S}_{\C}}
\newcommand{\Sp}{\mbox{Spin}}
\newcommand{\Spc}{\mbox{Spin}^{\C}}

\newtheorem{theorem}{Theorem}
\newtheorem{proposition}[theorem]{Proposition}
\newtheorem{lemma}[theorem]{Lemma}
\newtheorem{corollary}[theorem]{Corollary}

\newtheorem*{definition}{{\bf Definition}}

\begin{document}

\title{Survey on eigenvalues of the Dirac operator and geometric structures}
\author{Marcos Jardim  and Rafael F. Le\~ao \\ IMECC - UNICAMP \\
Departamento de Matem\'atica \\ Caixa Postal 6065 \\
13083-970 Campinas-SP, Brazil}

\maketitle

\begin{abstract}
We give a survey of results relating the restricted holonomy of a Riemannian spin manifold with lower bounds on the spectrum of its Dirac operator, giving a new proof of a result originally due to Kirchberg.
\end{abstract}

%\tableofcontents
%\newpage
%\baselineskip18pt

%---------------------------------------------------------------
%---------------------------------------------------------------

\section{Introduction} \label{intro}

Given a path connected, smooth Riemannian manifold $(M,g)$, its holonomy group is defined to be the
group of all linear transformations of the tangent space $T_pM$ induced by parallel
transport around loops based on $p\in M$. This Riemannian invariant encodes important
information about the manifold. In fact, Berger proved the following classification
result: if $(M,g)$ is neither locally a Riemannian product nor locally isometric to a
symmetric space, then its restricted holonomy group is one of the following:
\begin{enumerate}
\item $SO(n)$, with $\dim M=n$;
\item $U(n)$, with $\dim M=2n$;
\item $Sp(1)Sp(n)$, with $\dim M=4n\ge8$;
\item $SU(n)$, with $\dim M=2n\ge4$;
\item $Sp(n)$, with $\dim M=4n\ge8$;
\item $G_2$, with $\dim M = 7$;
\item $Spin(7)$, with $\dim M = 8$.
\end{enumerate}
The first case is the general one, while the other cases imply the existence of special
geometric structures on $M$. For instance, any manifold within case 3 is automatically Einstein, whilst
any manifold fitting in cases 4 through 7 is automatically Ricci-flat. For the proofs, we refer to \cite{J}.

Case 2 correspond to K\"ahler manifolds. Recall that a Riemannian manifold $(M,g)$ is K\"ahler
if it admits a complex structure, i.e. a bundle map $J:TM\to TM$ satisfying $J^2=-\id$, which is
parallel, i.e. $\nabla J=0$, where $\nabla$ denotes the Levi-Civitta connection.
The manifolds fitting in case 4 are Ricci-flat K\"ahler, also known in the literature
as Calabi-Yau manifolds. 

Case 3 correspond to quaternionic K\"ahler manifolds. Recall that a Riemannian manifold $(M,g)$ is 
quaternionic K\"ahler if ${\rm End}(TM)$ admits a parallel rank 3 sub-bundle $Q$ which is locally spanned
by almost complex structures $(I,J,K)$ satisfying the quaternionic relations, i.e. $IJ=K$, etc.
If the $(I,J,K)$ are globally defined and parallel, then $(M,g)$ is said to be a hyperk\"ahler manifold,
which are the manifolds of case 5.

Cases 6 and 7 are called exceptional cases, and compact examples were obtained by Joyce only recently
\cite{J}.

In this survey, we will show how the holonomy of a Riemannian spin manifold affects the spectrum of
the Dirac operator. More precisely, we will establish the following four theorems regarding the
smallest eigenvalue of the Dirac operator on manifolds of positive scalar curvature with various
holonomy groups.

In the general case, we have the following result due to Friedrich \cite{Fr}.

\begin{theorem}\label{t1}
Let $(M,g)$ be a compact Riemannian spin manifold of dimension $n$ and positive scalar curvature.
Then any eigenvalue $\lambda$ of the Dirac operator satisfies the following inequality:
$$ \lambda^2 \ge \frac{1}{4}\frac{n}{n-1}R_0 ~~, $$
where $R_0$ is the minimum of the scalar curvature. Moreover, if the equality is attained then
$(M,g)$ must be an Einstein manifold.
\end{theorem}

The case of K\"ahler manifolds was studied by Kirchberg, who proved the following theorem \cite{Ki1,Ki2}.

\begin{theorem}\label{t2}
Let $(M,J,g)$ be a compact K\"ahler spin manifold of real dimension $2n$ and positive scalar curvature.
Then any eigenvalue $\lambda$ of the Dirac operator satisfies the following inequality:
$$ \lambda^2 \ge \frac{1}{4}\frac{n+1}{n}R_0 ~~, $$
where $R_0$ is the minimum of the scalar curvature. Moreover, if the equality is attained then
$(M,g)$ must be K\"ahler-Einstein manifold and $n$ is odd.
\end{theorem}

The case of quaternionic K\"'ahler manifolds was considered by Kramer, Semmelmann and Weingart in
\cite{Kr1,Kr2}.

\begin{theorem}\label{t3}
Let $(M,I,J,K,g)$ be a compact quaternionic K\"ahler spin manifold of real dimension $4n$ and positive
scalar curvature. Then any eigenvalue $\lambda$ of the Dirac operator satisfies the following inequality:
$$ \lambda^2 \ge \frac{1}{4}\frac{(n+3)}{(4n+8)} R_0 ~~, $$
where $R_0$ is the minimum of the scalar curvature. Moreover, the equality is attained if and only if $M$
is a quaternionic projective space.
\end{theorem}

Finally, the situation for the remaining four cases is determined by the following result.

\begin{theorem}\label{t4}
The Dirac operator on a Riemannian spin manifold with restricted holonomy given by $SU(n)$, $Sp(n)$, $G_2$
or $Spin(7)$ has a nontrivial kernel.
\end{theorem}

The paper is organized as follows. We begin by reviewing some basic concepts and setting up notation in Section
\ref{sec_spin}. Theorem \ref{t1} is proved in Section \ref{f}, where we also comment on Theorem \ref{t4}. Section \ref{k} contains a new proof of Theorem \ref{t2}. We complete the paper with an overview of the proof of Theorem \ref{t3} in Section \ref{qk}.

\paragraph{\bf Acknowledgments.}
The first named author is partially supported by the CNPq grant number 300991/2004-5
and the FAPESP grant number 2005/04558-0. The second author's research was supported by
CNPq doctoral grant.

%---------------------------------------------------------------
%---------------------------------------------------------------

\section{Spin manifolds and the Dirac operator} \label{sec_spin}

Let $(M,g)$ be a smooth Riemannian manifold with dimension $n$. The bundle of orthonormal frames of $TM$ is a $SO_n$-principal bundle, and can be used to construct associated bundles. One bundle of particular interest that can be constructed in this way is the Clifford bundle. The standard action $\rho : GL(\R^n) \rightarrow R^n$ preserves the quadratic form of $R^n$ so this action can be naturally extended to an action $\rho : SO_n \rightarrow \Cl(\R^n)$.

\begin{definition}
The Clifford bundle is the vector bundle with standard fiber $\Cl_n$ given by
\begin{equation}
  \Cl(M) = P_{SO} \times_{\rho} \Cl(\R^n)
\end{equation}
where $\rho : SO_n \rightarrow \Cl(\R^n)$ is the action described above.
\end{definition}

Note that with this definition $\Cl(M)$ has a natural connection. In fact, we can look for the Levi-Civita connection of $(M,g)$ as a connection in the principal bundle $P_{SO_n}$, and, being $\Cl(M)$ an associated bundle, the connection on $P_{SO_n}$ induces a connection on $\Cl(M)$.

A bundle of modules for $\Cl(M)$ is a vector bundle $S$, with a Riemannian structure and a compatible connection, such that the fibers $S_p$ are modules over the fibers $\Cl(\R^n)_p$ of the Clifford bundle. 

\begin{definition}
Given a local orthonormal frame $\{ e_i \}$, we define de Dirac operator as the first order differential operator
\begin{equation}
  \begin{split}
    D : \Gamma( S ) \rightarrow \Gamma( S ) \\
    D \psi = \sum_{i=1}^n e_i \cdot \nabla_{e_i} \psi
  \end{split}
\end{equation}
\end{definition}

This definition works for every bundle of modules. In practice, however, we deal with bundles possessing further properties. It is natural to assume that the structures involved are compatible in some sense. With this in mind, we introduce the following definition.

\begin{definition}
A {\bf Dirac bundle} is bundle of modules $S$ with Riemannian structure and a compatible connection $\nabla$, such that
\begin{enumerate}
\item For vector fields $x,y \in \Cl(M)$, the Clifford action on $S$ is orthogonal, i.e.
\begin{equation}
  ( u \psi, u \phi ) = ( \psi , \phi ) ~~;
\end{equation} \\
\item The connection $\nabla$ on $S$ is a module derivation, i.e. for $s \in \Cl(M)$ and $\psi \in \Gamma(S)$ we have
\begin{equation}
    \nabla \left( s \cdot \psi \right) = \left( \nabla s \right) \cdot \psi + s \cdot \left( \nabla \psi \right)
\end{equation}
where $\nabla s$ denotes de connection of $\Cl(M)$ acting on $s$.
\end{enumerate}
\end{definition}

For certain Riemannian manifolds $(M,g)$ there exists a natural way to construct Dirac bundles. The relevant case for Dirac operators is the case of $\Sp$ manifolds. To understand what is a $\Sp$ manifold, let us look to the general case first.

\begin{definition}
Let $Q$ be a $SO_n$-principal bundle. A $\Sp$ structure on $Q$ is a $\Sp$-principal bundle $P$ and a double covering $\Lambda : P \rightarrow Q$ such that the diagram below be commutative
\begin{displaymath}
  \begin{CD}
    P \times \Sp_n @>>> P @>\pi>> M \\
    @VV\Lambda \times \lambda V  @VV\Lambda V @| \\
    Q \times SO_n @>>> Q @>\pi>> M
  \end{CD}
\end{displaymath}
where $\lambda : \Sp_n \rightarrow SO_n$ is the usual covering map.
\end{definition}

We say that a Riemannian manifold $(M,g)$ is $\Sp$ if the principal bundle of frames $P_{SO_n}$, associated with the tangent bundle, admits a $\Sp$ structure. Recall that $(M,g)$ admits a $\Sp$ structure if, and only if, the second Stiefel-Whitney class of its tangent bundle vanishes, $w_2(TM) = 0$.

\begin{definition}
The {\bf spinor} bundle of a $\Sp$ manifold $(M,g)$ is given by
\begin{equation}
  S = P_{\Sp} \otimes_{\rho} W
\end{equation}
where $W$ is a irreducible module for $\Cl(M)$
\end{definition}

The important fact is that the spinor bundle, as defined above, of a $\Sp$ manifold $(M,g)$ with the connection induced by the Levi-Civita connection is automatically a Dirac bundle. Thus, spinor bundles are a natural way to construct Dirac bundles. However, the explicit form of the action $\rho : \Cl(M) \rightarrow S$ is not clear in this construction. If the manifold $(M,g)$ is a complex manifold there is another way to look spinor bundles that make this action more evident.

In order to make this statement precise, we must understand complex $\Sp$ structures. If we consider the complexified algebra $\Cl(\R^n) \otimes \C$ we can look for $U(1)$ as being a subgroup of the units in this algebra. With this in mind, we define:

\begin{definition}
The $\Spc$ group is defined as the group
\begin{equation}
  \begin{split}
    \Spc = \Sp \times U(1) / \{ (-1, -1) \} \\
    \Spc \subset \Cl(\R^n) \otimes \C
  \end{split}
\end{equation}
\end{definition}

Using this group, we can define $\Spc$ structures in the following manner.

\begin{definition}
A $\Spc$ structure in a $SO_n$-principal bundle $Q$ is a $\Spc$-principal bundle $P$ and a covering map $\Lambda : P \rightarrow Q$ such that the diagram below is commutative
\begin{displaymath}
  \begin{CD}
    P \times \Spc_n @>>> P @>\pi>> M \\
    @VV\Lambda \times \lambda V  @VV\Lambda V @| \\
    Q \times SO_n @>>> Q @>\pi>> M
  \end{CD}
\end{displaymath}
\end{definition}

For every $\Spc$ structure there is an associated complex line bundle $\mc{L}$, often called the determinant of the $\Spc$ structure. The necessary topological condition for $(M,g)$ to admit a $\Spc$ structure is given in terms of this line bundle: an orientable manifold $(M,g)$ has a $\Spc$ structure if there exists a complex line bundle $\mc{L}$ such that
\begin{equation}
c^1(\mc{L}) =_{\mod 2} w_2(TM)
\end{equation}
where $c_1(\mc{L})$ denotes the first Chern class of $\mc{L}$. 

\begin{definition}
The bundle of complex spinors is defined as
\begin{equation}
  \Sc = P_{\Spc} \times_{\rho} W
\end{equation}
where $W$ is an irreducible module for $\Ccl_n = \Cl_n \otimes \C$, and $\rho : \Ccl_n \rightarrow W$ is the action induced by the inclusion $\Spc_n \subset \Ccl_n$.
\end{definition}

As in the case of spinors over a $\Sp$ structure, this bundle is a Dirac bundle; the relevant fact is that for complex manifolds we can give an explicit description of this bundle and of the action. Indeed, every complex manifold has a canonical $\Spc$ structure for which the determinant bundle is exactly the canonical bundle, $k_M$, of $M$; furthermore we have the identification
\begin{equation}
  \Sc \simeq \wedge^{0,*} M ~~.
\end{equation}
If we consider an unitary basis $\{ \xi^j , \bar{\xi}^j \}$ for $T^*M \otimes \C$, the action is explicit given by
\begin{equation}
  \begin{split}
    \rho(\xi^ j) \psi = - \sqrt{2} \xi^j \lrcorner \psi \\
    \rho(\bar{\xi}^j) \psi = \sqrt{2} \bar{\xi}^j \wedge \psi
  \end{split}
\end{equation}

For complex spin manifolds, we can construct both the spinor bundle $\mb{S}$ and the complex spinor bundle $\Sc$; they
are related in the following way.

\begin{proposition}
Let be $M$ a complex manifold with $\Sp$ structure. Let $\mb{S}$ be the spinor bundle associated to a given $\Sp$ structure and $\Sc$ the complex spinor bundle associated to the canonical $\Spc$ structure of $M$. Then
\begin{equation}\label{id}
  \Sc = \mb{S} \otimes k_M^{\frac{1}{2}}\simeq \wedge^{0,*}M \otimes k_M^{-\frac{1}{2}}
\end{equation}
where $k_M$ is the canonical bundle of $M$.
\end{proposition}

Now let us consider a hermitian vector bundle $E$ with connection $\nabla^A$. The bundles $\mb{S} \otimes E$, and $\Sc \otimes E$, have a natural module structure over $\Cl(M)$, defined simply in terms of the module structure of $\mb{S}$ or $\Sc$. Let $v \in \Cl(M)$ and $s \otimes t \in \mb{S} \otimes E$, then we have
\begin{equation}
  \rho(v) \left( s \otimes t \right) = \left( \rho(v) s \right) \otimes t
\end{equation}

It is easy to see that the bundle $\mb{S} \otimes E$ with the tensor product connection
$\nabla^{S \otimes A} = \nabla^S \otimes \mb{I} + \mb{I} \otimes \nabla^A$ is a Dirac bundle provided the connection
$\nabla^A$ is compatible with the hermitian structure of $E$. We can then define the {\bf twisted} Dirac operator as follows:
\begin{equation}
  \begin{split}
    D_A : \Gamma( \mb{S} \otimes E ) \rightarrow \Gamma( \mb{S} \otimes E ) \\
    D_A = \sum_i e^i \nabla^{S \otimes A}_i 
  \end{split}~~.
\end{equation}

The main tool usually employed in the study of the eigenvalues of Dirac operators is the Weitzenb\"ock formula. There are several variations of this formula depending on the case in question. For Dirac operators in a spinor bundle $\mb{S}$ associated to a $\Sp$ structure we have:
\begin{equation}\label{wf}
  D^2 = \Delta + \frac{1}{4} R ~~,
\end{equation}
where $R$ denotes the scalar curvature of $(M,g)$. If we consider the Dirac operator in a complex spinor bundle $\Sc$ associated to a $\Spc$ structure we have:
\begin{equation}
  D^2 = \Delta + \frac{1}{4} R + \frac{1}{2} F_{\sigma}
\end{equation}
where now $F_{\sigma}$ denotes the curvature 2-form of a fixed connection on the determinant bundle of the $\Spc$ structure. In the case of the canonical $\Spc$ structure this curvature is related to the curvature of $(M,g)$.

Now for twisted Dirac operators we must take into account the curvature of the connection $\nabla^A$ in $E$. If we consider $\mb{S} \otimes E$, where $\mb{S}$ is the spinor bundle associated to a $\Sp$ structure, we have:
\begin{equation}
  D_A^2 = \Delta^{S \otimes A} + \frac{1}{4} R + F_A ~~,
\end{equation}
where $F_A$ is the curvature 2-form of $\nabla^A$. Finally, if we consider $\Sc \otimes E$, where $\Sc$ is the complex spinor bundle associated to a $\Spc$ structure, then we have:
\begin{equation}
  D_A^2 = \Delta^{S \otimes A} + \frac{1}{4} R + \frac{1}{2} F_S + F_A ~~.
\end{equation}

%---------------------------------------------------------------
%---------------------------------------------------------------

\section{The Riemannian case}\label{f}

In this section, we show how to find sharp estimates for Dirac operators in Riemannian manifolds, a result first obtained by Friedrich. The idea is to consider a connection deformed using the module structure of Dirac bundles.

Let $E$ be a Dirac bundle over a Riemannian manifold $(M,g)$. In this bundle we can consider the deformed connection given by
\begin{equation}
  \nabla_v^f \psi = \nabla^A_v \psi + f v \cdot \psi
\end{equation}
Since $E$ is a Dirac bundle it is easy to see that the connection $\nabla^f$ is a metric connection, thus $E$ provided with the connection $\nabla^f$ still is a Dirac bundle.

To use this new connection to estimate eigenvalues of the Dirac operator we must find some kind of Weitzenb\"ock formula for the operators associated to $\nabla^f$. First, we define the deformed Dirac operator ($n=\dim M$):
\begin{equation}
  \begin{split}
    D^f &= \sum_i e_i \nabla^f_i = \sum_i e_i \nabla^A_i + \sum_i e_i B(e_i) \\
    &= D + f \sum_i e_i^2 = D - nf ~~.
  \end{split}
\end{equation}

The Laplacian associated to the connection $\nabla^A$ on $E$ is defined to be
\begin{equation}
  \Delta \psi = - \sum_i \nabla^A_i \nabla^A_i \psi - \sum_i div(e_i) \nabla_i \psi
\end{equation}
where $\nabla$ denotes the Levi-Civita connection of $(M,g)$ and $div(e_i)$ is given by
$div(e_i) = \sum_i g( \nabla_j e_i , e_j )$. If we consider an orthonormal basis $\{ e_i \}$
we have, using the compatibility of $\nabla$ with the metric, the following identity:
\begin{equation}\label{wdef}
\sum_i \nabla_i e_i = \sum_{ij} g( \nabla_i e_i , e_j ) e_j= -\sum_{ij} g( \nabla_i e_j , e_i ) e_j
= - \sum_j div(e_j) e_j ~~.
\end{equation}

\begin{lemma}
Let be $\Delta^f$ the Laplacian associated to the connection $\nabla^f$ and $\Delta$ the Laplacian of $\nabla^A$. Then
\begin{equation}
  \Delta^f = \Delta - 2 f D - \grad(f) + n f^2
\end{equation}
where $D$ is the Dirac operator on $E$ associated to $\nabla^A$
\end{lemma}
\begin{proof}
From the definition of the Laplacian we have
\begin{equation}
  \Delta^f \psi = - \sum_i \nabla^f_i \nabla^f_i \psi - \sum_i div(e_i) \nabla^f_i \psi ~~.
\end{equation}
The term $\sum_i \nabla^f_i \nabla^f_i \psi$ can be simplified:
\begin{equation}
  \begin{split}
    \sum_i \nabla^f_i \nabla^f_i \psi &= \sum_i ( \nabla^A_i + f e_i ) ( \nabla^A_i + f e_i ) \psi \\
    &= \sum_i \left( \nabla^A_i \nabla^A_i \psi + \nabla^A_i ( f e_i \psi ) + f e_i ( \nabla^A_i \psi ) - f^2 \psi
      \right) \\
    &= \sum_i \left( \nabla^A_i \nabla^A_i \psi + e_i ( \nabla_i f ) \psi + f ( \nabla_i e_i ) \psi + 2 f e_i ( 
      \nabla^A_i \psi ) - n f^2 \psi \right) \\
    &= \sum_i \nabla^A_i \nabla^A_i \psi + f \left( \sum_i \nabla_i e_i \right) \psi + \grad(f) \psi + 2 f D \psi - n 
      f^2 \psi \\
    &= \sum_i \nabla^A_i \nabla^A_i \psi - f \left( \sum_j div(e_j) e_j \right) \psi + \grad(f) \psi + 2 f D \psi - n 
      f^2 \psi
  \end{split}
\end{equation}

On the other hand, we can write $\sum_i div(e_i) \nabla^f_i \psi$ as:
\begin{equation}
  \sum_i div(e_i) \nabla^f_i \psi = \sum_i div(e_i) \nabla^A_i \psi + f \sum_i div(e_i) e_i \psi ~~.
\end{equation}
Equation (\ref{wdef}) now follows easily from the last three equations.
\end{proof}

\begin{lemma} \label{ide_def}
For the deformed connection $\nabla^f$, we have the Weitzenb\"ock formula
\begin{equation}\label{ide_def_eq}
  \left( D - f \right)^2 = \Delta^f + \mc{F} + (1-n) f^2 ~~,
\end{equation}
where $\mc{F}$ is curvature 2-form of the connection on the Dirac bundle in question.
\end{lemma}
\begin{proof}
First, note that if $f$ is a function on $M$ then $D( f \psi ) = \grad(f) \psi + f D \psi$,
since $\nabla$ is a derivation and $\grad(f) = \sum_i e_i \nabla^A_i f = \sum_i e_i e_i(f)$.
It then follows that
\begin{equation}
  \left( D - f \right)^2 = D^2 - 2fD - \grad(f) + f^2 \label{d-f} ~~.
\end{equation}

Combining equations (\ref{lap_def}) and (\ref{d-f}) we obtain
\begin{equation}
  \left( D - f \right)^2 = \Delta^f + \left( D^2 - \Delta \right) + (1-n) f^2 ~~,
\end{equation}
thus (\ref{ide_def_eq}) follows from the application of the usual
Weitzenb\"ock formula to this last equation.
\end{proof}

We are finally ready to prove our first main result, Theorem \ref{t1}. Take $\psi$ such that $D \psi = \lambda \psi$. Making the deformation parameter $f$ constant and equal to $\frac{\lambda}{n}$, the equation (\ref{ide_rie}) takes
the form
\begin{equation}
  \lambda^2 \left( \frac{n-1}{n} \right) \psi = \Delta^{\frac{\lambda}{n}} \psi + \frac{1}{4} R \psi
\end{equation}
Now take inner product with $\psi$ to obtain
\begin{equation}
  \lambda^2 \left( \frac{n-1}{n} \right) \mid \mid \psi \mid \mid^2_{L_2} = \mid \mid \nabla^{\frac{\lambda}{n}} \psi
    \mid \mid^2_{L_2} + \frac{1}{4} \int_M R \mid \psi \mid^2
\end{equation}
Since $\mid \mid \nabla^{\frac{\lambda}{n}} \psi \mid \mid^2_{L_2} \geq 0$ and estimating $R \geq R_0$, we can conclude that
\begin{equation}\label{ineq}
  \lambda^2 \geq \frac{1}{4} \frac{n}{n-1} R_0
\end{equation}
as desired.

\begin{proposition}
If exists a section $\psi \in \mb{S}$, such that
$$ D \psi = \frac{1}{4} \frac{n}{n-1} R_0 \psi $$
then the scalar curvature $R$ is constant and we have
\begin{equation} \label{eq_kil_res}
  \nabla_x \psi = \mp \frac{1}{2} \sqrt{\frac{R_0}{n(n-1)}} x \psi
\end{equation}
for any $x \in TM$. 
\end{proposition}

\begin{proof}
In order for the equality in (\ref{ineq}) to hold, we must have $R = R_0$ and
\begin{equation}
\nabla^{\frac{\lambda}{n}} \psi = 0 ~~,
\end{equation}
and the Proposition follows easily.
\end{proof}

Motivated by the previous Proposition, we introduce the following Definition.

\begin{definition}
A Killing spinor $\psi$ is a spinor that satisfies the equation
\begin{equation} \label{eq_kil}
  \nabla_x \psi = \mu x \cdot \psi ~~ \forall x \in TM
\end{equation}
for some constant $\mu$.
\end{definition}

Riemannian manifolds admitting a Killing spinor have strong geometrical properties, see \cite[Section 5.2]{F}.

\begin{proposition}
Let $(M,g)$ be a Riemannian manifold with $\Sp$ structure and let $\psi$ be a Killing spinor. Then $(M,g)$ is an Einstein manifold and we have $\mu^2 = \frac{1}{4} \frac{1}{n(n-1)} R$. Moreover, if $\mu\neq0$, then $(M,g)$ is locally irreducible and has constant sectional curvature.
\end{proposition}

Therefore, it is immediate to conclude that if the lower bour in (\ref{ineq}) is actually attained, then $M$ must be an Einstein manifold, as desired.
(...)

%---------------------------------------------------------------
%---------------------------------------------------------------

\section{The K\"ahler case}\label{k}

The estimate for the general Riemannian case obtained in the previous section can't be satisfied for K\"ahler manifolds. If some section $\psi \in \mb{S}$ satisfies the equation $\nabla_x \psi = \mu x \cdot \psi$ is easy to see that this section satisfies

\begin{equation}
  D \psi = -2 n \mu \psi
\end{equation}
where $n$ is the real dimension of $(M,g)$.

But if $M$ is a K\"ahler manifold we can use the K\"ahler form to construct another eigensection of D, using the section $\psi$. In other words, if $\psi$ is an eigensection of $D$, with eigenvalue $\lambda$, then the section $\omega \psi$ is another eigensection of $D$. But the eigenvalue associated to this section is $\lambda^{\prime}= \frac{2n-4}{2n} \lambda$.

This immediately implies that in a K\"ahler manifold with real dimension different from $2$, we can't have a spinor satisfying the equality in the Friedrich estimate.

To obtain a sharp estimate we must modified the deformation introduced by Friedrich. Let $(M,g,J)$ be a K\"ahler manifold with $\Sp$ structure. Let $\mb{S}$ be the spinor bundle associated to this $\Sp$ structure. We know that the Levi-Civita connection of $(M,g,J)$ induces a unique connection in $\mb{S}$. In terms of this connection we define the deformed connection
\begin{equation}
  \nabla^{a,b}_x \psi = \nabla_x \psi + a x \cdot \psi + i b J(x) \cdot \alpha( \psi )
\end{equation}

The first term $a x \cdot \psi$ is exactly the Friedrich deformation. The second term $b J(x) \cdot \alpha( \psi )$ involves the complex structure of $M$ and the parity operator on $\mb{S}$. To understand the parity operator remember that $M$, besides the $\Sp$ structure, also have a canonical $\Spc$ structure. So we can consider the spinor bundle $\Sc$ associated to this $\Spc$ structure. The two spinor bundles are related by

\begin{equation}
  \Sc = \mb{S} \otimes k_M^{\frac{1}{2}}
\end{equation}
where $k_M$ is the canonical bundle of $M$. So we can write the spinor bundle $\mb{S}$ as

\begin{equation}
  \mb{S} = \Sc \otimes k_M^{-\frac{1}{2}} \simeq \wedge^{0,*}M \otimes k_M^{-\frac{1}{2}}
\end{equation}

Now the parity operator on forms, $\alpha$, is given by $\alpha( \psi_p ) = (-1)^p \psi_p$ for $\psi_p \in \wedge^{0,p}M$, and using the above description for $\mb{S}$ we immediately seen that this operator is well defined on $\mb{S}$.

The parity operator $\alpha$ on $\mb{S}$ can be related to the K\"ahler structure of $M$ in a suitable way. The K\"ahler form defines a splitting of $\mb{S}$ that naturally defines operators related to the parity operator. To see how this happens remember that the action of $\Cl(M)$ on $\Sc \simeq \wedge^{0,*}M$ is given by
\begin{equation}\label{eq_ac_sp}
\rho(\xi^ j) \psi = - \sqrt{2} \xi^j \lrcorner \psi  ~~{rm and}~~
\rho(\bar{\xi}^j) \psi = \sqrt{2} \bar{\xi}^j \wedge \psi
\end{equation}
can be extended to $\mb{S}$ in a natural way.

Using this action and the fact that the K\"ahler form can be written as

\begin{equation}
  \omega = i \sum_{i=1}^n \xi^i \wedge \bar{\xi}^i
\end{equation}
we immediately have

\begin{proposition}
Let $\omega$ be the K\"ahler form of $M$ see as an operator on $\mb{S}$. Let $\psi_p \in \wedge^{0,p}M \otimes k^{-\frac{1}{2}}$, then we have
\begin{equation}
  \omega \psi_p = i(2p-n) \psi_p
\end{equation}

So we can write $\mb{S}$ as a sum of eigenbundles of $\omega$

\begin{equation}
  \mb{S} = \oplus_p \mb{S}_p
\end{equation}
where $\mb{S}_p = \wedge^{0,p}M \otimes k^{-\frac{1}{2}}$
\end{proposition}

Using this decomposition of $\mb{S}$ we can define a square root for the parity operator of $\mb{S}$. Now with the decomposition $\mb{S} = \oplus_p \mb{S}_p$ we define

\begin{equation}
  \mc{I} = \sum_{k=0}^{n} (i)^k p_k
\end{equation}
and note that $\mc{I}^2 = \alpha$.

On spinors $\psi \in \mb{S}$, the K\"ahler form and the complex structure $J$ of $M$ are related
by the following two Lemmata.

\begin{lemma} \label{lem_j_kah}
Let $\alpha$ be a 1-form. Then we have the relation

\begin{equation}
  \alpha \omega - \omega \alpha = 2 J( \alpha )
\end{equation}
\end{lemma}

\begin{proof}
Being $\alpha$ a 1-form we have the identity

\begin{equation}
  \alpha \omega - \omega \alpha = - 2 \alpha \lrcorner \omega
\end{equation}

In other way

\begin{equation}
  \begin{split}
    \alpha \lrcorner \omega (y) &= \omega( \alpha^{\flat} , y ) = g( \alpha^{\flat} , J( y ) ) \\
    &= - g( J( \alpha^{\flat} ) , y ) = -J( \alpha ) (y)
  \end{split}
\end{equation}
\end{proof}

\begin{lemma} \label{lem_for_kah}
The K\"ahler form, see as an operator on $\mb{S}$, satisfies the relation
\begin{equation}
  \sum_{i=1}^n J( e^i ) e^i = 2 \omega
\end{equation}
\end{lemma}

\begin{proof}
Using the previous Lemma we have
\begin{equation}
\begin{split}
\sum_{i=1}^n j( e^i ) e^i &= \frac{1}{2} \sum_{i=1}^n ( e^i \omega - \omega e^i ) e^i
= \frac{1}{2} \left[ n \omega + \sum_{i=1}^n e^i \omega e^i \right] \\
    &= \frac{1}{2} \left[ (4-n) \omega + n \omega \right] = 2 \omega
  \end{split}
\end{equation}
\end{proof}

To compute the Laplacian associated to the connection $\nabla^{a,b}$ we need to introduce the deformed Dirac operator

\begin{equation}
  \tilde{D} = \sum_{i=1}^{2n} J(e^i) \nabla_i
\end{equation}

It is interesting to note that this operator and the Dirac operator $D$ are related by the operator $\mc{I}$. It is easy to see that

\begin{equation}
  \tilde{D} = - \mc{I} D \mc{I}^* = \mc{I}^* D \mc{I}
\end{equation}
where $\mc{I}$ is the formal adjoint of $\mc{I}$. Besides, the following relations hold:
\begin{equation}\label{rels}
D \mc{I} D \mc{I} = \mc{I} D \mc{I} D  ~~,~~ \mc{I} D \mc{I} = - \mc{I}^* D \mc{I}^*
\end{equation}

Using all this relations we are able to compute the Laplacian associated to the deformed connection $\nabla^{a,b}$

\begin{theorem} \label{teo_lap_kir}
The Laplacian associated to the connection $\nabla^{a,b}$ is given by

\begin{equation}
  \Delta^{a,b} \psi = \Delta \psi +n(a^2 + b^2) \psi - 2 a D \psi - 2i b \tilde{D} \alpha( \psi ) + 4iab \omega 
    \alpha( \psi )
\end{equation}
\end{theorem}

The proof consist in write the Laplacian and manipulate the expression using the above identities. This is a huge calculation without any insights and will be omitted.

To use this Laplacian in estimates for the eigenvalues of the Dirac operator, we need to control the terms $2i b \tilde{D}$ and $4iab \omega$.

%Using that $\mc{I}^2 = \alpha$ we immediately see that $D \mc{I}^2 = - \mc{I}^2 D$. This fact, with the relation of $D$ and $\tilde{D}$, allow us to write the relations (\ref{rels})
%\begin{equation}
%  \begin{split}
%    D \mc{I} D \mc{I} = \mc{I} D \mc{I} D \\
%    \mc{I} D \mc{I} = - \mc{I}^* D \mc{I}^*
%  \end{split}
%\end{equation}

Let $E_{\lambda}(D)$ denotes the eigenspace of $D$ with eigenvalue $\lambda$.

\begin{proposition}
If $\psi \in E_{\lambda}(D)$ the expression

\begin{equation}
  e_{\lambda} \psi = ( D + \lambda ) \mc{I}^* \psi
\end{equation}
define an endomorphism $e_{\lambda}: E_{\lambda}(D) \rightarrow E_{\lambda}(D)$ such that

\begin{equation}
  e_{\lambda}^4 + 4 \lambda^4 = 0
\end{equation}
\end{proposition}

\begin{proof}
Taking $\psi \in E_{\lambda}(D)$ and using that $D^2$ commutes with $\mc{I}^*$ we have

\begin{equation}
\begin{split}
D( e_{\lambda} \psi ) &= D^2 \mc{I}^* \psi + \lambda D \mc{I}^* \psi 
= \lambda^2 \mc{I}^* \psi + \lambda D \mc{I}^\psi \\
&= \lambda \left( \lambda + D \right) \mc{I}^* \psi
\end{split}\end{equation}

So $e_{\lambda}$ really defines an endomorphism of $E_{\lambda}(D)$. In other way, supposing that $\psi \in E_{\lambda}(D)$, and using the above identities, we have

\begin{equation}
  \begin{split}
    e_{\lambda}^2 \psi &= ( D \mc{I}^* + \lambda \mc{I}^* ) ( D \mc{I}^* + \lambda \mc{I}^* ) \psi \\
    &= D \mc{I}^* D \mc{I}^* \psi + \lambda D (\mc{I}^*)^2 \psi + \lambda \mc{I}^* D \mc{I}^* \psi + \lambda^2 
      (\mc{I}^*)^2 \psi \\
    &= - D \mc{I} D \mc{I} \psi - \lambda (\mc{I}^*)^2 D \psi - \lambda \mc{I} D \mc{I} \psi + \lambda (\mc{I}^*)^2 
      \psi \\
    &= - D \mc{I} D \mc{I} \psi - \lambda \mc{I} D \mc{I} \psi \\
    &= -\mc{I} D \mc{I} D \psi - \lambda \mc{I} D \mc{I} \psi \\
    &= -2 \lambda \mc{I} D \mc{I} \psi
  \end{split}
\end{equation}

Repeating the same calculation to $e_{\lambda}^3 \psi$ and $e_{\lambda}^4 \psi$ we have

\begin{equation}
e_{\lambda}^3 \psi = - 2 \lambda^2 ( D + \lambda ) \mc{I} \psi ~~ {\rm and} ~~
e_{\lambda}^4 \psi = -4 \lambda^4 \psi ~~. 
\end{equation} \end{proof}

In particular, the expression $e_{\lambda}^4 = - 4 \lambda^4$ says that the only possible eigenvalues of $e_{\lambda}$ on $E_{\lambda}(D)$ are the complex numbers $\pm ( 1 \pm i ) \lambda$. With this in mind we define

\begin{equation}
  E_{\lambda}^k(D) = \{ \psi \in E_{\lambda}(D) \mid e_{\lambda} \psi = i^k ( 1 + i ) \lambda \psi \}
\end{equation}

\begin{corollary} \label{cor_e_lam}
Let $\lambda$ be an eigenvalue of $D$. Then exists some $k \in \{ 0,1,2,3 \}$ and $\psi \in E_{\lambda}^k(D)$, with $\psi \neq 0$. Beside this, $\psi$ satisfies

\begin{equation}
  \tilde{D} \psi = - \lambda \left( i^k (1+i) \mc{I} - 1 \right) \psi
\end{equation}
\end{corollary}

\begin{proposition} \label{prop_rel_pro}
Let $\lambda \neq 0$ be an eigenvalue of $D$, and let $\psi \in E_{\lambda}^k(D)$. Then the projection operators relative to the decomposition  $\mb{S} = \oplus_j \mb{S}_j$ satisfies

\begin{equation}
  \begin{split}
    \mid \mid p_{4 l -k -1} \psi \mid \mid = \mid \mid p_{4l -k} \psi \mid \mid \\
    p_{4l-k+1} \psi = p_{4l-k+2} \psi = 0
  \end{split}
\end{equation}
\end{proposition}

\begin{proof}
Using the explicit action in terms of $\{ \xi^i , \bar{\xi}^i \}$ we have

\begin{equation}
  p_j J( x ) - J( x ) p_{j-1} = i ( p_j x - x p_{j-1} )
\end{equation}

This implies that

\begin{equation}
  p_j \tilde{D} - \tilde{D} p_{j-1} = i ( p_j D - D p_{j-1} )
\end{equation}

Using the fact that $\tilde{D}$ is self-adjoint, we have, for $\psi \in E_{\lambda}^k(D)$, that

\begin{equation}
  \begin{split}
    \langle p_j \tilde{D} \psi \mid \psi \rangle - \langle \psi \mid p_{j-1} \tilde{D} \psi \rangle = \\
    i \langle (p_j D - D p_{j-1}) \psi \mid \psi \rangle = \\
    i \langle p_j D \psi \mid \psi \rangle + i \langle p_{j-1} \psi \mid D \psi \rangle = \\
    i \lambda \langle p_j \psi \mid \psi \rangle + i \lambda \langle p_{j-1} \psi \mid \psi \rangle = \\
    i \lambda \left( \mid \mid p_j \psi \mid \mid^2 - \mid \mid p_{j-1} \psi \mid \mid^2 \right)
  \end{split}
\end{equation}

Using this and the corollary (\ref{cor_e_lam}) we finally get

\begin{equation}
  (i + i^{k+j}) \mid \mid p_j \psi \mid \mid^2 = \left( i+(-i)^{k+j} \right) \mid \mid p_{j-1} \psi \mid \mid^2
\end{equation}

Now the result follows.
\end{proof}

This relations allow us to control the terms involving $\tilde{D}$ and $\omega$ in the Laplacian.

\begin{proposition} \label{prop_con}
Let $\lambda \neq 0$ be an eigenvalue of $D$, and $\psi \in E_{\lambda}^k(D)$. Then we have

\begin{equation}
  \begin{split}
    \langle - i \tilde{D} \psi \mid \mc{I}^2 \psi \rangle = (-1)^{k+1} \lambda \mid \mid \psi \mid \mid^2 \\
    (-i \omega \psi \mid \mc{I}^2 \psi \rangle = (-1)^k \mid \mid \psi \mid \mid^2
  \end{split}
\end{equation}
\end{proposition}

\begin{proof}
We know that $\mc{I}^2 = \alpha$ and that $D \alpha = -\alpha D$. Then it is immediate that $\langle \psi \mid \mc{I}^2 \psi \rangle = 0$. With this we have

\begin{equation}
  \langle -i \tilde{D} \psi \mid \mc{I}^2 \psi \rangle = - i^k(1-i) \lambda \langle \psi \mid \mc{I} \psi \rangle
\end{equation}

Using proposition (\ref{prop_rel_pro}) we have

\begin{equation}
  \begin{split}
    \langle \psi \mid \mc{I} \psi \rangle &= \sum_j (-i)^j \mid \mid p_j \psi \mid \mid^2 \\
    &= \sum_j (-i)^{4j-k-1} \mid \mid p_{4j-k-1} \psi \mid \mid^2 + \sum_j (-i)^{4j-k} \mid \mid p_{4j-k} \psi \mid 
      \mid^2 \\
      &= i^k (1+i) \sum_j \mid \mid p_{4j-k} \psi \mid \mid^2 = \frac{1}{2} i^k(1+i) \mid \mid \psi \mid \mid^2
  \end{split}
\end{equation}

The last two equations gives

\begin{equation}
  \langle - i \tilde{D} \psi \mid \mc{I}^2 \psi \rangle = (-1)^{k+1} \lambda \mid \mid \psi \mid \mid^2
\end{equation}

Now for $\omega$ write

\begin{equation}
  \omega = i \sum_j (2j-n) p_j
\end{equation}
and this implies that

\begin{equation}
  \begin{split}
    \langle - i \omega \psi \mid \mc{I}^2 \psi \rangle &= \langle \sum_j (2j-n) p_j \psi \mid \sum_k (-1)^k p_k \psi 
      \rangle \\
      &= \sum_j (-1)^j (2j-n) \mid \mid p_j \psi \mid \mid^2 \\
      &= \sum_p (-1)^{4p-k}(8p-2k-n) \mid \mid p_{4p-k} \psi \mid \mid^2 \\
      &+ \sum_p (-1)^{4p-k-1}(8p-2k-2-n) \mid \mid p_{4p-k-1} \psi \mid \mid^2 \\
      &= 2 (-1)^k \sum_p \mid \mid p_{4p-k} \psi \mid \mid^2 = (-1)^k \mid \mid \psi \mid \mid^2
  \end{split}
\end{equation}
\end{proof}

\begin{theorem}
Let $M$ be a K\"ahler manifold with $\Sp$ structure and let $D$ be the associated Dirac operator. Then, if $\lambda$ is an eigenvalue of $D$, $\lambda$ satisfies

\begin{equation}
  \lambda^2 \geq \frac{1}{4} \frac{n+2}{n} R_0
\end{equation}
where $R_0$ denotes the minimum of the scalar curvature of $M$.
\end{theorem}

\begin{proof}
This is an immediate consequence of the above considerations if we take $a = \frac{\lambda}{n+2}$ and $b = (-1)^{k+1} \frac{\lambda}{n+2}$.
\end{proof}

If the equality is satisfied we can prove, in the same way that was proved for the Riemannian case, that the manifold $(M,g)$ is an Einstein manifold with constant scalar curvature. But in the K\"ahler is another consequences, if the equality is satisfied, using properties of the projection operators we can prove that the manifold $M$ must have odd complex dimension.

A more general argument involving twistor operators, that obtain a sharp estimate for the case of even complex dimension was found by Kirchberg in \cite{Ki2}.

%------------------------------------------------------
%------------------------------------------------------

\section{The quaternionic K\"ahler case}\label{qk}

In this section we only will give the idea of the proof, which can be found in \cite{Kr1}. As in the K\"ahler case, the idea is to consider further structures of the manifold $M$ to obtain a better estimate. In the K\"ahler case, the K\"ahler structure was considered in terms of the decomposition of the spinor bundle in eigenbundles of the K\"ahler form, and using this we were able to deform the connection and the respective Weitzenb\"ock formula to obtain a sharp estimate.

In the quaternionic K\"ahler case, the idea is similar. Kraines \cite{Kr} proved that for quaternionic K\"ahler manifolds there exists a fundamental 4-form $\Omega$, which can be used to decompose the spinor bundle \cite{Hi}. This decomposition can then be used to obtain a sharp estimate, but for quaternionic K\"ahler manifolds there exists an alternative argument that leads to the same decomposition.

As we mentioned in the Introduction, quaternionic K\"ahler manifolds are characterized by having holonomy group $Sp(1) Sp(n)$. Now representation theory for $Sp(1) Sp(n)$ can be used. As it is well know, c.f. \cite{Sa}, all representations of $Sp(1) Sp(n)$ are given in terms of the fundamental representation $H = \mb{H} \simeq \C^2$ and $E = \mb{H}^n \simeq \C^{2n}$.

Let $M$ a quaternionic K\"ahler manifold, and let $\Hh$ and $\E$ be the vector bundles associated to the fundamental representations defined above. The fact that all representations of $Sp(1) Sp(n)$ can be given in terms of the fundamental representations implies that all vector bundles with structure group $Sp(1) Sp(n)$ can be given in terms of the vector bundles $\Hh$ and $\E$. In particular one can prove, c.f. \cite{Sa}, that the complexified tangent bundle of $M$ can be written as
\begin{equation}
  TM = \Hh \otimes \E ~~,
\end{equation}
while the spinor bundle is given by, c.f. \cite{Kr1}:
\begin{equation}
  \mb{S} = \oplus_{r=0}^n Sym^r \Hh \otimes \wedge_0^{n-r} \E ~~,
\end{equation}
where $\wedge_0^{n-r} \E$ denotes some subspace of $\wedge^{n-r} \E$ determined by the action of $Sp(1) Sp(n)$.

These decompositions of the tangent bundle of $M$ and of the spinor bundle are equivalent to the decompositions of the tangent bundle and spinor bundle of a K\"ahler manifold in terms of eigenbundles of the complex structure and eigenbundles of the K\"ahler form $\omega$. In fact, $\mb{S}_r = Sym^r \Hh \otimes \wedge_0^{n-r} \E$ are precisely
the eigenbundles of the fundamental 4-form of $M$. Besides, the Clifford multiplication can be described similarly to
(\ref{eq_ac_sp}).

With these descriptions, a Weitzenb\"ock formula adapted to quaternionic K\"ahler manifolds can be derived, and Theorem \ref{t3} is a direct consequence of this formula just as in the proofs of the previous results for Riemannian and for K\"ahler manifolds.

%------------------------------------------------------
%------------------------------------------------------

\end{document}